\documentclass[a4paper,
      abstracton]{scrartcl}
\pdfoutput=1

\usepackage[tracking, kerning, spacing]{microtype}
\microtypecontext{spacing=nonfrench}

\newcommand{\institute}[1]{}
\usepackage{amsmath,amsthm,amssymb}
\newtheorem{theorem}{Theorem}
\newtheorem{conjecture}{Conjecture}
\theoremstyle{definition}
\newtheorem{definition}{Definition}
\theoremstyle{remark}
\newtheorem{example}{Example}
\newtheorem{remark}{Remark}

\usepackage{amstext,amsfonts,amssymb,amscd,amsbsy,amsmath}
\usepackage{tikz}

\newcommand{\RR}{\mathbb{R}}
\newcommand{\ZZ}{\mathbb{Z}}
\newcommand{\poly}{\Pi}

\newcommand{\define}[1]{\emph{#1}}

\usepackage{xspace}
\usepackage{url}

\begin{document}

\newcommand{\polymake}{\emph{polymake}\xspace}

\title{Cellular sheaf cohomology in Polymake}
\author{Lars Kastner, 
Kristin Shaw, 
Anna-Lena Winz
}
\institute{Lars Kastner \at Freie Universit\"at Berlin, Arnimallee 3, 14195 Berlin, Germany, \email{k.l@fu-berlin.de}
\and Kristin Shaw \at Technische Universit\"at Berlin, Institut f\"ur Mathematik, Stra\ss e des 17. Juni 136, 10623 Berlin, Germany, \email{shaw@math.tu-berlin.de}
\and Anna-Lena Winz \at Freie Universit\"at Berlin, Arnimallee 3, 14195 Berlin, Germany, \email{anna-lena.winz@fu-berlin.de}
}

\maketitle
\titlehead

\newcommand{\weirdabstract}[1]{\abstract*{#1}\abstract{#1}}

\begin{abstract}
This chapter provides a guide to our \polymake extension  \textbf{cellularSheaves}. 
We first define cellular sheaves on polyhedral complexes in Euclidean space, as well as cosheaves, and their (co)homologies. 
As motivation, we summarise some  results from toric and tropical geometry linking cellular sheaf cohomologies to 
cohomologies of algebraic varieties. 
We then give an overview of the structure of the extension
\textbf{cellularSheaves} for \polymake. Finally, we illustrate the usage of the
extension with  examples from toric and tropical geometry.
\end{abstract}

\section{Introduction}

Given a polyhedral complex $\poly$ in $\mathbb{R}^n$,  a cellular sheaf (of
vector spaces) on $\poly$ associates to every face of $\poly$ a vector space
and to every face relation a map of the associated vector spaces (see
Definition \ref{def:sheaf}).  Just as with usual sheaves, we can compute the
cohomology of  cellular sheaves (see Definition \ref{def:homologies}).  The
advantage is that cellular sheaf cohomology  is the cohomology of a  chain
complex consisting of finite dimensional vector spaces (see Definition
\ref{def:cochain}).

Polyhedral complexes often arise  in algebraic geometry. Moreover, in some
cases, invariants of algebraic varieties can be recovered as cohomology of
certain cellular sheaves on polyhedral complexes.  Our two major classes of
examples come from toric and tropical geometry and are presented in Section
\ref{sec:thms}.  We refer the reader to \cite{Fulton:toric} for a guide to
toric geometry and polytopes.  For an introduction to tropical geometry see
\cite{BIMS}, \cite{dianeBernd}.

The main motivation for this polymake extension was to 
 implement tropical homology, as introduced by
Itenberg, Katzarkov, Mikhalkin, and Zharkov in \cite{trophom}. Tropical
homology is the homology of particular cosheaves which can be defined on 
any   polyhedral complex. 
When the polyhedral complex arises as the tropicalisation of a family of complex projective
varieties, the tropical homology groups give information about the Hodge numbers 
of a generic member of the family. 
However, this is just one particular instance of cellular (co)sheaf (co)homology which our 
package is capable of dealing with. 
Cellular (co)sheaves have also been a powerful tool in recent years in the field of  
applied topology, notably in persistent homology, sensor networks and network coding  \cite{Ghrist}, \cite{Curry}.  With this  \polymake extension, 
(co)sheaves on polyhedral complexes can be constructed from scratch. We hope that this will allow for a range of 
uses of the extension beyond just the ones from  combinatorial algebraic geometry that we point out here. 

The definitions of cellular (co)sheaves and their (co)homologies are given in Section \ref{sec:defs}. 
Then a description of our implementation of  cellular sheaves and their cohomologies
in \polymake is given in Section \ref{sec:implementation}.  The framework
\polymake already provides a large number of combinatorial objects and tools,
making it easy to construct the polyhedral complexes of interest. In 
Section \ref{sec:code}, we  illustrate the usage of the extension  in a
variety of examples from tropical and algebraic geometry. Finally, Section \ref{sec:future} outlines
some potential future directions and applications for our extension.

\section{Cellular sheaf cohomology}\label{sec:defs}

We begin by defining cellular sheaves and cosheaves, as well as their
cohomologies and homologies.  We provide an explicit example of a sheaf and a
cosheaf  which have been implemented in our package.  In Section \ref{sec:thms}
we give an overview of some theorems relating these sheaf cohomologies and
cosheaf homologies to the cohomology of some algebraic varieties. 

\subsection{Definitions}

A \define{polyhedral complex} $\poly$ is a finite collection of polyhedra in
$\mathbb{R}^n$ with the property that any face of a polyhedron in $\poly$ is
also in $\poly$ and the intersection of any two polyhedra in $\poly$ is a face
of both. 

Let $\poly^{i}$ denote the collection of polyhedra in $\poly$ of dimension $i$.
For polyhedra $\sigma, \tau \in \poly$, we use the notation $\tau \leq \sigma$
to denote that \define{$\tau $ is a face of $\sigma$}.

\begin{definition}\label{def:orient}
Given a polyhedral complex $\poly$ and a chosen orientation of each polyhedron in $\poly$ define the \define{orientation map} for each $i$, 
$$\mathcal{O} : \poly^{i-1} \times \poly^{i} \to \{0, -1, +1\}$$
by 
\[
    \mathcal{O}(\tau, \sigma)\ :=\  
\begin{cases}
    0 & \text{if } \tau \not \leq \sigma \\
    +1 & \begin{minipage}{8cm}if the orientation of $\tau \subset \partial \sigma$ coincides with that of $\tau$\end{minipage} \\           
    -1 &  \begin{minipage}{8cm}if the orientation of $\tau \subset \partial \sigma$ differs from that of $\tau$.\end{minipage} \\         
\end{cases}
\]
\end{definition}

A polyhedral complex $\poly$ can be considered as a category where the objects
are the polyhedra and the morphisms are given by inclusions. i.e. $(f : \tau
\to \sigma) \in \text{Mor}(\poly)$ if and only if $\tau \leq \sigma$. We use the
notation $\poly^{op}$ to denote the category obtained from $\poly$ by using the
same objects and reversing the directions of all morphisms. 

We will be interested in cellular sheaves and cosheaves of vector spaces on polyhedral
complexes. Viewing $\poly$ as a category as described above, we can give a
succinct definition of cellular (co)sheaves.
Let $\text{Vect}_K$ denote the category of vector spaces over a field $K$.

\begin{definition}\label{def:sheaf}
Let $\poly$ be a polyhedral complex, then a \define{cellular sheaf} $\mathcal{G}$ and a \define{cellular cosheaf} $\mathcal{F}$ are functors 
$$ \mathcal{G}: \poly \to \text{Vect}_K \quad \quad \quad \quad \mathcal{F}: \poly^{op} \to \text{Vect}_K.$$ 
\end{definition}

To summarise, this means that a cellular sheaf consists of the following data:

\begin{itemize}
\item for each polyhedron \(\sigma\) in $\poly$ a vector space $\mathcal{G}(\sigma)$ and,
\item given $\tau, \sigma \in
\poly$ satisfying $\tau \leq \sigma$, a morphism $\rho_{\tau \sigma}: \mathcal{G}(\tau) \to
\mathcal{G}(\sigma)$.
\end{itemize}
In particular, for $\gamma\le\tau\le\sigma$ the restriction morphisms commute, i.e. we have
\[
\rho_{\gamma\sigma}\ =\ \rho_{\tau\sigma}\circ\rho_{\gamma\tau}.
\]

A cellular cosheaf is similar except that the morphisms
are in the opposite direction $\iota_{\sigma \tau}:  \mathcal{F}(\sigma) \to \mathcal{F}(\tau)$. 

A sheaf of vector spaces in the usual sense is a contravariant functor from the
category of open sets of a topological space to $\text{Vect}_K$ which satisfies
additional axioms.  A polyhedral complex can be equipped with a finite topology
known as the Alexandrov topology, and the above definition of cellular sheaf as
a functor produces a sheaf in the usual sense in this topology.  Due to the
simplicity of the cellular sheaves and the Alexandrov topology, no additional
sheaf axioms are required. The reader is directed to \cite[Chapter 4]{Curry}
for more details. 

\begin{example}\label{ex:constantsheaf}
As a first example we can define the constant sheaf by setting \(\mathcal{G}(\sigma)\) to be the one dimensional vector space $K$
and the maps $\rho_{\tau \sigma}: \mathcal{G}(\tau) \to \mathcal{G}(\sigma)$ to be the identity for all 
\(\tau, \sigma \in \poly\) such that $\tau \leq \sigma$. 
The constant cosheaf can be defined in a similar fashion. 
\end{example}

\begin{example}\label{ex:Wsheaves}
Let \(\poly\) be a polyhedral complex in \(\RR^n\). For \(\sigma \in \poly\) set \(L(\sigma)\)
to be the linear subspace of \(\RR^n\) parallel to the face \(\sigma\). 
For \( p \in \ZZ_{\ge 0}\) we define the sheaf \(W^p\) as follows:
\[
W^p(\sigma) = \wedge^p \ L(\sigma)
\]
and the maps $\rho_{\tau \sigma}: W^p(\tau) \to W^p(\sigma)$ are given by the wedges of the natural inclusions \(L(\tau) \to L(\sigma)\) for $\tau \leq \sigma$.
By convention, the  sheaf \(W^0\) is the constant sheaf from Example \ref{ex:constantsheaf}.
\end{example}

\begin{example}\label{ex:Fcosheaves}
Next we give an example of a cosheaf on a polyhedral complex.
The homology of this particular cosheaf is the tropical homology from \cite{trophom}, and will come up in many examples in 
Sections \ref{sec:Berg} and \ref{sec:hypersurface}. 

Let \(\poly\) be a polyhedral complex in \(\RR^n\), then we define
\[F_p(\sigma) = \sum_{\sigma < \gamma} \wedge^p L(\gamma) .\]
If \(\tau \le \sigma\), then \(\{\gamma \ | \ \sigma < \gamma\} \subset \{ \gamma \ |\  \tau < \gamma\}\), so we get 
a natural inclusion \( \iota_{\sigma \tau}: F_p(\sigma) \to F_p(\tau)\).
Analogously to Example \ref{ex:Wsheaves}, we obtain the constant cosheaf from Example \ref{ex:constantsheaf} for $p=0$. 

\end{example}

Note that by dualising the vector spaces $\mathcal{G}(\sigma)$ for all $\sigma$ we can transform a cellular sheaf $\mathcal{G}$ into
a cellular cosheaf in a natural way. A sheaf can be created from a cosheaf in the analogous way.

For a given (co)sheaf we build (co)chain complexes in the following two parallel definitions. The definitions originally appeared in this 
form in 
\cite{Curry}.

\begin{definition}\label{def:cochain}
Given a polyhedral complex $\poly$ and a cellular sheaf $\mathcal{G}$, define
the \define{cellular cochain groups} and \define{cellular cochain groups with
compact support}, respectively,  as 

$$C^{q}(\poly ;\  \mathcal{G}) := \bigoplus_{\substack{ \dim \sigma = q  \\ \sigma \text{ compact}}} \mathcal{G}(\sigma) \quad \text{and} \quad C_{c}^{q}(\poly ;\  \mathcal{G}) := \bigoplus_{\dim \sigma = q} \mathcal{G}(\sigma).$$
The cellular cochain maps (usual or  with compact support) $$d: C^{q}(\poly ;\  \mathcal{G}) \to 
C^{q+1}(\poly ;\  \mathcal{G}) \quad \text{and} \quad d: C^{q}_{c}(\poly ;\  \mathcal{G}) \to 
C^{q+1}_{c}(\poly ;\  \mathcal{G})$$ 
are given component-wise
for \(\tau \in \poly^{q}\) and \(\sigma \in \poly^{q+1}\) 
by  $d_{\tau\sigma}:\mathcal{G}(\tau)\to\mathcal{G}(\sigma)$, where
\[
d_{\tau\sigma}\ :=\ \begin{cases}
\mathcal{O}(\tau, \sigma)\cdot\rho_{\tau\sigma} & \tau\le\sigma\\
0 & \mbox{else}.
\end{cases}
\]
 \end{definition}

\begin{definition}\label{def:chain}

Given a polyhedral complex $\poly$ and a cellular cosheaf $\mathcal{F}$ define
the \define{cellular chain groups} and the \define{Borel-Moore cellular chain
groups}, respectively, as

$$C_{q}(\poly ;\  \mathcal{F}) := \bigoplus_{\substack{\dim \sigma = q \\ \sigma \text{compact}}} \mathcal{F}(\sigma)
\quad \text{and} \quad 
C^{BM}_{q}(\poly ;\  \mathcal{F}) := \bigoplus_{\dim \sigma = q} \mathcal{F}(\sigma).$$
The cellular chain maps (usual or Borel-Moore) $$\partial : C_{q}(\poly
;\  \mathcal{F}) \to C_{q-1}(\poly ;\  \mathcal{F}) \quad \text{and} \quad \partial : C_{q}^{BM}(\poly
;\  \mathcal{F}) \to C_{q-1}^{BM}(\poly ;\  \mathcal{F})$$
are given component-wise 
for \(\sigma \in \poly^{q}\) and \(\tau \in \poly^{q-1}\) 
by $\partial_{\sigma\tau}:\mathcal{F}(\sigma)\to\mathcal{F}(\tau)$, where
\[
\partial_{\sigma\tau}\ :=\ \begin{cases}
\mathcal{O}(\sigma, \tau)\cdot\iota_{\sigma \tau} & \sigma\ge\tau\\
0 & \mbox{else}.
\end{cases}
\]
 
\end{definition}

\begin{remark}
It may seem counter intuitive that the usual cellular cochains are supported
only on compact faces and cellular cochains with compact support are supported
on all faces.
After we define cohomology of a cellular sheaf and homology of a cellular
cosheaf below, a sanity check can be performed to compute the cohomology of the
constant sheaf from Example \ref{ex:constantsheaf} on your favourite
non-compact polyhedral complex. Under some reasonable conditions on the
polyhedral complex, the cellular cohomology of the constant sheaf will be
isomorphic to the ordinary singular cohomology of the polyhedral complex. The
analogous statement holds for the compactly supported versions. See
\cite[Example 6.2.4]{Curry} for a simple example and~more details. 
\end{remark}

\begin{definition}\label{def:homologies}
The cellular sheaf cohomology (with compact support) of $\mathcal{G}$ is the
cohomology of the cellular cochain complex (with compact support) from
Definition \ref{def:cochain}.

The cellular (Borel-Moore) cosheaf homology of $\mathcal{F}$ is the homology of
the cellular (Borel-Moore) chain complex from Definition \ref{def:chain}. 
\end{definition}

\subsection{Connections with classical algebraic geometry}\label{sec:thms} 
 
In this section we present some particular connections between cellular
(co)ho\-mo\-lo\-gies of certain (co)sheaves and cohomology of complex algebraic
varieties.  Explicit demonstrations of the statements of these theorems are
given along with the \polymake code in Section \ref{sec:code}.
 
To a rational polytope $\Delta \subset \mathbb{R}^n$ we can associate a toric
variety $TV(\Delta)$ by first considering its normal fan then building the
toric variety corresponding to this fan following \cite{Fulton:toric}.  The
following theorem relates the cohomology of the sheaves of $p$-differential
forms $\Omega^p$ on $TV(\Delta)$ with the cohomologies of the sheaves $W^p$
from Example \ref{ex:Wsheaves} on the polytope $\Delta$.

\begin{theorem}\label{danilov:differentials}\cite[Remark 12.4.1]{Danilov} 
Let $\Delta \subset \mathbb{R}^n$ be a rational polytope and $TV(\Delta)$ its
associated toric variety. Then 
\[
H^q( TV(\Delta);\  \Omega^p ) \cong H^q(\Delta;\  W^p) \otimes_{\mathbb{R}}
\mathbb{C}.
\]
\end{theorem}
 
In particular, when the rational polytope $\Delta$ is simple, the associated
toric variety is smooth and $H^q(TV(\Delta);\  \Omega^p)  \cong H^{p,
q}(TV(\Delta))$, where the last vector space denotes the $(p, q)^{\text{th}}$
part of the Hodge decomposition of the cohomology of $TV(\Delta)$. 
As an example in Section \ref{sec:polytope},  we compute $H^q(\Delta;\  W^p)$ for $\Delta$ a three dimensional 
cube using the extension \textbf{cellularSheaves}.

The next two connections relate to the $F$-cosheaves defined in Example
\ref{ex:Fcosheaves}. We begin with an arrangement of hyperplanes $\mathcal{A}$
in $\mathbb{C}P^d$. It is a theorem of Orlik and Solomon that the cohomology
of the complement of the arrangement depends only on the combinatorics of the
arrangement, in other words the corresponding matroid (see 
\cite{Katz:Matroid} or \cite{Oxley} for an introduction to matroid theory). Moreover, for any
matroid there is a combinatorially described Orlik-Solomon algebra, which
provides the cohomology ring of the complement when we are in the situation
above (see \cite{Orlik}). 

To any matroid $M$ we can associate a fan $B(M)$ in Euclidean space known as
the Bergman fan of $M$ (see \cite{Ardila:Berg}). If $M$ is the matroid of an
arrangement $\mathcal{A}$ the fan $B(M)$ is the tropicalisation of the
complement $\mathbb{C}P^d \backslash \mathcal{A}$ under a suitable embedding to
a complex torus. 

\begin{theorem}[\cite{Zharkov}]\label{thm:berg}
Let $M$ be any matroid, $B(M)$ its Bergman fan, and $F_p$ the cosheaves from
Example \ref{ex:Fcosheaves} on $B(M)$, then
$$OS^p(M) = F_p(v)^*,$$ 
where $OS^p$ denotes the $p^{\text{th}}$ graded part of the Orlik-Solomon
algebra of $M$, and $v$ is the vertex of the fan. 
\end{theorem}

In Section \ref{sec:Berg}, we illustrate how we can compute the dimensions of
the graded pieces of the Orlik-Solomon algebra of a matroid using the \polymake
applications \texttt{matroid}, \texttt{tropical} and the extension
\textbf{cellularSheaves}.

Lastly, the statements relating the cohomology of the complements of
arrangements and the tropical homology of matroidal fans in Theorem
\ref{thm:berg} can be generalised and even refined in the setting of
tropicalisation of projective complex algebraic varieties. We state the theorem
and refer the reader to \cite{trophom} for the precise definitions of smooth
$\mathbb{Q}$-tropical projective varieties and tropical limits.

\begin{theorem}\cite{trophom}
Consider a $1$-parameter family of complex projective varieties $\pi: \mathcal{X} \to D^o$, where $D^o$ is the punctured disc. 
Suppose that the tropical limit $\text{Trop}(\mathcal{X}) = X$ is a smooth $\mathbb{Q}$-tropical projective variety. Then 
$$\dim H^{p, q}(\mathcal{X}_t) = \dim H_{q}(X;\  F_p),$$
where $H^{p, q}(\mathcal{X}_t) $ is the $(p, q)^{\text{th}}$ part of the Hodge decomposition of 
$\mathcal{X}_t = \pi^{-1}(t)$ for a generic $t \in D^o$. 
\end{theorem}

It is due to the above theorem that the homology of the $F$-cosheaves is also known as tropical homology. 

As of yet, \polymake does not have compact tropical varieties as objects.
Therefore, in the examples presented in Section \ref{sec:code}, we do not
produce Hodge numbers of complex projective algebraic varieties, but rather
Betti numbers of limit mixed Hodge structures of non-complete varieties. Examples of
this can be found in \cite{DanKho86} for hypersurfaces and complete
intersections, also in \cite{KatzStapledon} from a more tropical point of view.
Future plans to implement tropical homology for compact (and hence also
projective) tropical varieties are outlined in Section \ref{sec:compact}. 

\section{Implementation in \polymake}
\label{sec:implementation}

The framework \polymake is a mathematical software for polyhedral geometry.
Its objects of interest are mainly combinatorial, such as cones, polyhedra,
graphs, fans and polyhedral complexes. Toric and tropical geometry provide many
ways for using \polymake to solve computational tasks from algebraic
geometry. In particular, \polymake provides the applications \texttt{tropical}
for tropical geometry and \texttt{fulton} for toric geometry. Furthermore,
\polymake interfaces several other software packages which may be useful in
our context, such as \emph{gfan} (\cite{gfan}) for tropical computations and
\emph{Singular} (\cite{Singular}) for algebraic geometry.  See \cite{polymake:tropical} 
for an overview of 
the most current implemented \polymake features for  tropical geometry. 

The interface language of \polymake is perl. For improved performance one can
write and attach C++ code. The combinatorial objects are realised as objects
with properties, e.g. the object \texttt{Polytope} has the properties
\texttt{VERTICES} and \texttt{F\_VECTOR} amongst many others. Since solving
certain problems can be very expensive time- and resource-wise, \polymake
adheres to the principle of lazy evaluation: properties are only computed when
needed and then stored with the object, so they do not have to be recomputed.

Computation of properties is done via \polymake's internal rule structure. A
rule takes a certain set of input properties and then computes a certain set of
output properties. 
When asked for a certain property of an object, \polymake creates a queue of
rules to apply in order to get this property from any set of given properties,
if this is possible.

Take for example the following code snippet:

\begin{verbatim}
object PolyhedralFan {

   property ORIENTATIONS : Map<Set<Set<Int> >, Int>;

   rule ORIENTATIONS: HASSE_DIAGRAM, FAN_DIM, RAYS,
      LINEALITY_SPACE{
		... # Code
   }

}
\end{verbatim}

Here the object \texttt{PolyhedralFan} is equipped with a new property
\texttt{ORIEN\-TA\-TIONS} which one needs for computing tropical homology, see
Definition \ref{def:orient}. Then a rule is created, that computes \texttt{ORIENTATIONS}
from the properties \texttt{HASSE\_DIA\-GRAM}, \texttt{FAN\_DIM}, \texttt{RAYS}
and \texttt{LINEALITY\_SPACE} of the \texttt{PolyhedralFan}.

Internally in \polymake, every polyhedral complex $\Pi$ in $\mathbb{R}^n$ is considered as a polyhedral fan $\Sigma$ in $ 
\mathbb{R}^{n+1}$ 
intersected with the hyperplane defined by  $x_0 = 1$.  
Every  face of $\Pi$ is indexed by a subset of the rays of $\Sigma$. 
The one dimensional faces of $\Sigma$ whose direction $\vec{v}=(v_0,\ldots,v_n)$ satisfies $v_0 = 1$
correspond to vertices of $\Pi$. The one dimensional faces of $\Sigma$ whose direction $\vec{v}$ satisfies $v_0 = 0$
correspond to unbounded  one dimensional faces of $\Pi$. 

\begin{definition}
Let $\Sigma$ be a fan in $\mathbb{R}^{n+1}$  and $\Pi = \Sigma \cap \{x_0 = 1\}$. 
A far vertex of $\Pi$ is a ray of $\Sigma$ whose direction $\vec{v}$ satisfies $v_0 = 0$. 
 A face $\sigma$ of $\Pi$ is, 
\begin{itemize}
\item  a far face if its index set consists only of far vertices;
\item a non far face if its index set contains at least one non-far vertex; 
\item a  bounded face  if it contains no far vertices;
 \item an unbounded face if it is neither a  far face nor a bounded face. 
 \end{itemize}
 
Our extension \textbf{cellularSheaves} adds the properties \texttt{FAR\_FACES}, \texttt{BOUNDED\_FACES}, \texttt{UNBOUNDED\_FACES} to a polyhedral complex. 
\end{definition}

Computing 
orientations for the polyhedral fan avoids complications caused by the different types of faces of the polyhedral complex. 
 Since the object 
\texttt{Po\-ly\-he\-dral\-Com\-plex} is derived from the object \texttt{PolyhedralFan},
 it will have the property \texttt{ORIENTATIONS} as well.

\subsection{Obtaining the \textbf{cellularSheaves} extension}
The extension can be installed on a Linux system with the most recent
\polymake version with the following two steps. First clone the repository
with
\begin{verbatim}
git clone \
   http://www.github.com/lkastner/cellularSheaves \
   FOLDER
\end{verbatim}
into a folder named \texttt{FOLDER}. Second start \polymake, and import the extension using
\begin{verbatim}
import_extension("FOLDER");
\end{verbatim}

The extension introduces the new objects \texttt{Sheaf} and \texttt{CoSheaf} from 
Definition \ref{def:sheaf}. See Section \ref{sec:sheavesImp} for more details on the implementation. 
A basic usage scenario looks like
\begin{verbatim}
application "fan";
$pc = new PolyhedralComplex( 
   check_fan_objects(new Cone(cube(3))));
$w1 = $pc->wsheaf(1);
\end{verbatim}

First we switch to the application \texttt{fan}, since this is the application
our extension adds functionality to. The next line takes the three dimensional
cube and turns it into a polyhedral complex. Then we ask for the $W^1$-sheaf of
Example \ref{ex:Wsheaves}.

We implemented most methods dealing with pure linear algebra in C++. The file 
\begin{verbatim}
apps/fan/include/linalg.h
\end{verbatim}
contains the C++ code.  These linear algebra methods, especially those
assembling a chain complex from given block matrices, perform significantly
better when implemented in C++ than in perl.

\subsection{Sheaves and cosheaves}\label{sec:sheavesImp}

In our extension we introduce the objects \texttt{Sheaf} and \texttt{CoSheaf}.
As implemented, these objects  have two properties. The first is a map  from a
collection of sets of integers to matrices. This property represents the vector
spaces of a (co)sheaf.  The second is a map from pairs of
sets in this collection to matrices. These matrices represent the morphisms
between these vector spaces. 

The vector spaces and morphisms are stored in the following two properties of a
(co)sheaf:
\begin{verbatim}
property BASES : Map<Set<Int>, Matrix>;

property BLOCKS : Map<Set<Set<Int> >, Matrix >;
\end{verbatim}

Let us rephrase this in terms of the Definitions \ref{def:cochain} and
\ref{def:chain}. Let $\poly$ be a polyhedral complex with a sheaf
$\mathcal{G}$, and let $\tau\le\sigma$ be a face relation in $\poly$. 
The faces of $\Pi$ are encoded as index sets of the rays of vertices of the defining polyhedral 
fan $\Sigma$. 
As an object
in \polymake, the sheaf $\mathcal{G}$ has the property \texttt{BASES} containing
the bases of the vector spaces $\mathcal{G}(\gamma)$ for all $\gamma$ in $\Pi$. 
Also
for the sheaf $\mathcal{G}$, the property \texttt{BLOCKS} contains a
matrix representing the map $\rho_{\tau \sigma}: \mathcal{G}(\tau)\to\mathcal{G}(\sigma)$, for each pair of faces $\tau \leq
\sigma$. This
matrix is written using the bases from the property \texttt{BASES}.
Analogously for cosheaves, the property \texttt{BASES} of $\mathcal{F}$
contains the bases of $\mathcal{F}(\tau)$ for all $\tau \in \Pi$. The property
\texttt{BLOCKS} contains a matrix representing the map
$\iota_{\sigma \tau}: \mathcal{F}(\sigma )\to \mathcal{F}(\tau)$, for each pair of faces $\tau \leq
\sigma$.

For the purpose of computing sheaf cohomology it is also necessary to store
morphisms for certain non-face relations, these will then consist of zero
matrices of the appropriate sizes. 

The main (co)sheaf constructors are \texttt{fcosheaf}, \texttt{wsheaf} which 
produce the (cosheaves) from   Examples \ref{ex:Wsheaves} and \ref{ex:Fcosheaves}. 
 These are user methods attached to a polyhedral complex.
Each of these methods takes a non-negative integer as parameter, which determines
the $p$ in the wedge power for the $F$-cosheaves and $W$-sheaves.

\subsection{Chain complexes and homologies}

The last new important objects are chain complexes introduced as
\texttt{ChainComplex}. A chain complex comes with the properties
\texttt{DIFFERENTIALS}, \texttt{BETTI\_NUMBERS}, \texttt{HOMOLOGY} and
\texttt{IS\_WELLDEFINED}. It can be created by giving an array of matrices as
the property \texttt{INPUT\_DIFFERENTIALS}. Furthermore it has a user method
\texttt{print()} providing a human readable sequence format of the chain
complex.
Dually, we introduce the object
\texttt{CoChainComplex}. In reality this is just a wrapper around the
object \texttt{ChainComplex} for the user's convenience. 

Currently there are two (co)homology methods in our extension for a given (co)sheaf. 
They differ by which faces are considered when building the chain
complex. 
\begin{enumerate}
\item \texttt{usual\_chain\_complex}: This method considers only the bounded
faces of the given polyhedral complex, i.e. it computes $C_{\bullet}(\poly ;\  \mathcal{F})$.
\item \texttt{borel\_moore\_complex}: This method uses all non-far faces of a
given polyhedral complex, i.e. it computes $C^{BM}_{\bullet}(\poly ;\  \mathcal{F})$.

Analogously
\item \texttt{usual\_cochain\_complex} gives $C^{\bullet}(\poly ;\ \mathcal{G})$ and
\item \texttt{compact\_support\_complex} gives $C^{\bullet}_{c}(\poly ;\ \mathcal{G})$.
\end{enumerate}

\section{Examples and usage}\label{sec:code}
This section provides  sample code and output for some specific examples. 
These examples are chosen so as to highlight the connections to cohomology of complex 
algebraic varieties described in Section \ref{sec:thms}.

\subsection{Polytopes}\label{sec:polytope}
We consider the polyhedral complex that consists of a three dimensional cube $C$ and 
all its faces. We will compute the $W$-sheaves for $C$ as well as the Betti numbers of the 
cohomology groups $H^q(C ;\  W^p)$ for all $p, q$ from $0$ to $3$. 
\begin{verbatim}
application "fan";
$pc = new PolyhedralComplex( 
   check_fan_objects(new Cone(cube(3))));
@betti = ();
for(my $i=0; $i<4; $i++){
   my $w = $pc->wsheaf($i);
   my $s = $pc->usual_cochain_complex($w);
   push @betti, $s->BETTI_NUMBERS;   
}
print new Matrix(@betti);
\end{verbatim}

\noindent The first step turns the three dimensional cube into a polyhedral complex. Then we loop over all possible $W$-sheaves and save the Betti numbers in a matrix. This results in

\begin{verbatim}
fan > print new Matrix(@betti);
1 0 0 0
0 3 0 0
0 0 3 0
0 0 0 1
\end{verbatim}
We see that \(\dim H^q(C;\  W^p) = 0\) if \(p \not= q\). 
The diagonal \(\dim H^p(C;\  W^p)\) is the dual h-vector of the polytope that defines the polyhedral complex.
This relationship holds for any simple polytope $\Delta$. See  \cite[Corollary, pg.~6]{Brion} for the statement in terms 
of the dual fan of $\Delta$. 
\begin{verbatim}
fan > $cube = polytope::cube(3);

fan > print $cube->DUAL_H_VECTOR;
1 3 3 1
\end{verbatim}
The toric variety of the three dimensional cube is $ X = \mathbb{P}^1 \times \mathbb{P}^1  \times \mathbb{P}^1$. Notice that  $\text{dim} H^{p} ( X; \Omega^p ) = 1 $ if $p = 0, 3$, $\text{dim} H^{p} ( X ; \ \Omega^p) = 3$ if $p = 1, 2$ and $\text{dim} H^{q} ( X ; \ \Omega^p)  = 0$ for $p\neq q$. 

\subsection{Bergman fans and tropical linear spaces}\label{sec:Berg}

We can build the Bergman fan $B(M)$ of a matroid $M$ and compute
the usual homology of the $F$-cosheaf - this means that we 
only consider bounded faces.
Assuming the matroid is connected, the Bergman fan of a matroid has 
a unique bounded face which is the vertex $v$. 
Therefore, the cellular chain groups $C_{q}(B(M);\   F_p)$ are $ 0$ unless $q = 0$. 

In the following examples we will see that
$$\dim H_{0}(B(M);\  F_p) = \dim OS^p(M),$$
where $OS^p(M)$ is the $p^{th}$ graded part of the Orlik-Solomon algebra of
$M$. This follows from Theorem \ref{thm:berg}. 
Notice that when  $v$ is  the vertex of the Bergman fan,  
$$H_{0}(B(M);\  F_p) = C_{0}(B(M);\  F_p) = F_p(v)$$
and $H_{i}(B(M);\  F_p) = 0$ for $i >0$. 

When $M$ is a rank $d+1$ matroid on $n+1$ elements arising from a non-central hyperplane
arrangement $\mathcal{A}$ in $\mathbb{C}P^{d}$, the Orlik-Solomon algebra is
isomorphic to the cohomology ring of the complement $\mathcal{C} :=
\mathbb{C}P^{d} \backslash \mathcal{A}$ of the arrangement. There is a
canonical embedding of $\mathcal{C} \to(\mathbb{C}^*)^{n}$, and the
tropicalisation of this is the Bergman fan of the matroid. Therefore, we see
that the homology of the $F$-cosheaf on a tropicalisation recovers
cohomological information about the original variety. 

\begin{example}\label{ex:tropicalLine}
Our first example is to compute the tropical homology of a tropical line in
$\mathbb{R}^2$. This is the tropicalisation of a generic line $L \subset
\mathbb{C}^2$ intersected with the torus $(\mathbb{C}^*)^2$. Notice that this
space is homeomorphic to $\mathbb{C} P^1 \backslash \{p_1, p_2, p_3\}$, so 
that $\dim H^{0}(L \cap (\mathbb{C}^*)^2 ;\  \mathbb{C}) = 1$ and 
$\dim H^{1}(L \cap (\mathbb{C}^*)^2 ;\  \mathbb{C}) = 2$. 
The
tropical line is the Bergman fan of the uniform matroid of rank $2$ on $3$
elements. 

We start by creating this polyhedral complex in \polymake:
\begin{verbatim}
application "matroid";
$m = uniform_matroid(2,3);
application "tropical";
$t = matroid_fan<Max>($m);
$t->VERTICES;
application "fan";
$berg = new PolyhedralComplex($t);
\end{verbatim}
Next, we construct the associated $F$-cosheaves up to the  dimension of the Bergman fan and compute their usual chain complexes:
\begin{verbatim}
$f0 = $berg->fcosheaf(0);
$f1 = $berg->fcosheaf(1);
$s0 = $berg->usual_chain_complex($f0);
$s1 = $berg->usual_chain_complex($f1);
\end{verbatim}
Now we ask for the Betti numbers and obtain:
\begin{verbatim}
fan > print $s0->BETTI_NUMBERS;
1 0
fan > print $s1->BETTI_NUMBERS;
2 0
\end{verbatim}
We can also compute the Borel-Moore homology. Here every face of the Bergman fan 
contributes to the Borel-Moore chain groups (see Definition \ref{def:homologies}).

\begin{verbatim}
$bm0 = $berg->borel_moore_complex($f0);
$bm1 = $berg->borel_moore_complex($f1);
\end{verbatim}
gives
\begin{verbatim}
fan > print $bm0->BETTI_NUMBERS;
0 2
fan > print $bm1->BETTI_NUMBERS;
0 1
\end{verbatim}
Notice that we obtain 
\(\dim H_q(B(M);\  F_p)= \dim H_{d-q}^{BM}(B(M);\  F_{d-p})\) for $d = 1$, which is the dimension of the Bergman fan.
This is the homological version of Poincar\'e duality for matroidal fans and
tropical manifolds from \cite{JSS}. 

\end{example}

\begin{example}\label{ex:K4}
In this example we will study the Bergman fan of the matroid of the
complete graph on four vertices. This is the matroid of the so-called braid
arrangement of lines in $\mathbb{C}P^2$, whose complement is the moduli space
of $5$-marked genus $0$ curves $\mathcal{M}_{0, 5}$ (see \cite{Ardila:Berg}). We use the applications ``graph" and 
``matroid" to first construct the Bergman fan. 

\begin{verbatim}
application "graph";
$g = complete(4);
application "matroid";
$m = matroid_from_graph($g);
application "tropical";
$t = matroid_fan<Max>($m);
$t->VERTICES;
application "fan";
$berg = new PolyhedralComplex($t);
\end{verbatim}
We compute the usual and the Borel-Moore homology of the 
$F$-cosheaf.
\begin{verbatim}
@betti_usual = ();
@betti_bm = ();
for(my $i=0; $i<3; $i++){
   my $f = $berg->fcosheaf($i);
   my $s = $berg->usual_chain_complex($f);
   my $bm = $berg->borel_moore_complex($f);
   push @betti_usual, $s->BETTI_NUMBERS;
   push @betti_bm, $bm->BETTI_NUMBERS;
}
\end{verbatim}
This gives the following Betti numbers:
\begin{verbatim}
fan > print new Matrix(@betti_usual);
1 0 0
5 0 0
6 0 0

fan > print new Matrix(@betti_bm);
0 0 6
0 0 5
0 0 1
\end{verbatim}
Again we see that we have \(\dim H_q(B(M); \ F_p)= \dim H_{d-q}^{BM}(B(M);
\ F_{d-p})\), where now $d=2$. 

\end{example}

\begin{example}\label{ex:troplinearspace}
A tropical linear space is not necessarily a fan. Nevertheless the Betti
numbers of the tropical homology of the tropical linear space and of its
recession fan agree. In this example, we start with the Bergman fan of the
uniform matroid of rank $3$ on $6$ elements and compare its homology with that
of the tropical linear space of a valuated matroid with the aforementioned
matroid as its underlying matroid.

\begin{verbatim}
$m = matroid::uniform_matroid(3,6);
$t = tropical::matroid_fan<Max>($m);
$t->VERTICES;
application "fan";
$berg = new PolyhedralComplex($t);
@betti_usual = ();
@betti_bm = ();
for(my $i=0; $i<3; $i++){
   my $f = $berg->fcosheaf($i);
   my $s = $berg->usual_chain_complex($f);
   my $bm = $berg->borel_moore_complex($f);
   push @betti_usual, $s->BETTI_NUMBERS;
   push @betti_bm, $bm->BETTI_NUMBERS;
}
\end{verbatim}
gives
\begin{verbatim}
fan > print new Matrix(@betti_usual);
1 0 0
5 0 0
10 0 0

fan > print new Matrix(@betti_bm);
0 0 10
0 0 5
0 0 1
\end{verbatim}
Now we consider a valuated matroid whose underlying matroid is uniform of rank $3$
on $6$ elements and construct the corresponding
tropical linear space. 

\begin{verbatim}
$v = [0,0,3,1,2,1,0,1,0,2,2,0,3,0,4,1,2,2,0,0];
$val_matroid = new matroid::ValuatedMatroid<Min>( 
   BASES=>matroid::uniform_matroid(3,6)->BASES, 
   VALUATION_ON_BASES=>$v,N_ELEMENTS=>6);
$tls = tropical::linear_space($val_matroid);
@betti_usual = ();
@betti_bm = ();
for(my $i=0;$i<3;$i++){
   my $fi = $tls->fcosheaf($i);
   my $si=$tls->usual_chain_complex($fi);
   my $bmi=$tls->borel_moore_complex($fi);
   push @betti_usual, $si->BETTI_NUMBERS;
   push @betti_bm, $bmi->BETTI_NUMBERS;
}
\end{verbatim}

\noindent returns, 

\begin{verbatim}
fan > print new Matrix(@betti_usual);
1 0 0
5 0 0
10 0 0

fan > print new Matrix(@betti_bm);
0 0 10
0 0 5
0 0 1
\end{verbatim}

\noindent which is the same as for the Bergman fan of the matroid above. 

\end{example}

\begin{example}\label{ex:WsheavesLinSpace}
This example demonstrates that the usual cohomology and the compactly supported 
cohomology of the $W$-sheaves on tropical linear spaces seem to satisfy some interesting and also
potentially useful vanishing theorems. We demonstrate this with a single
example and then summarise the vanishing phenomena in Conjecture
\ref{conj:vanishing} below.  We continue with the same tropical linear space
from Example \ref{ex:troplinearspace}. 

\begin{verbatim}
@wbetti_usual = ();
@wbetti_cs = ();
for(my $i=0;$i<3;$i++){
   my $wi = $tls->wsheaf($i);
   my $wsi=$tls->usual_cochain_complex($wi);
   my $wcsi=$tls->compact_support_complex($wi);
   push @wbetti_usual, $wsi->BETTI_NUMBERS;
   push @wbetti_cs, $wcsi->BETTI_NUMBERS;
}
\end{verbatim}

\noindent returns, 

\begin{verbatim}
fan > print new Matrix(@wbetti_usual);
1 0 0
0 4 0
0 0 1

fan > print new Matrix(@wbetti_cs);
0 0 10
0 0 32
0 0 28

\end{verbatim}

\end{example}
\begin{conjecture}\label{conj:vanishing}
Let $L \subset \mathbb{R}^n$ be a tropical linear space of dimension $d$. Then we have 
$$H^{q}(L;\ {W}^p) = 0 \quad \text{if } p \neq q \quad \quad \text{ and } \quad \quad H^{q}_c(L;\ {W}^p) = 0 \quad \text{if } q \neq d.$$ 

To date,  the above conjecture has been checked on all realisable tropical linear spaces in 
$\text{Trop}(\text{Gr}(3, 6))$  using our package. 

By considering the Euler characteristics of the complexes $C^{\bullet}(L; \ {W}^p)$  and $C_c^{\bullet}(L;\ {W}^p)$  we have: 

$$(-1)^pH^{p}(L\ ; \ {W}^p) = \sum_{ q = 0}^d (-1)^q\binom{q}{p}f^b_q,$$

$$(-1)^d H^{d}_c(L;\ {W}^p) = \sum_{ q = 0}^d (-1)^q\binom{q}{p}f_q$$
where $f^b = (f_0^b, \dots f_d^b)$ is the $f$-vector of the bounded faces of
$L$ and $f = (f_0, \dots f_d)$ is the $f$-vector of $L$.  If the above conjecture holds, then understanding 
the 
$f$-vector of a tropical linear space comes down to understanding the possible
dimensions of $H^{q}_{\bullet}(L ;\ {W}^p)$. For example, it is possible to
bound the $f^b$-vector by bounding $H^{p}(L;\ {W}^p)$. This would give an
approach to the $f$-vector conjecture for tropical linear spaces (see \cite{Speyer}) similar to the proof of the upper bound conjecture for polytopes. 
 
\end{conjecture}

\subsection{Tropical hypersurfaces}\label{sec:hypersurface}

Using the a-tint package (\cite{Hampe}) we can construct tropical hypersurfaces
in \polymake from piecewise  integer affine functions which are convex,
i.e.~tropical polynomials. 
These examples demonstrate how one can start directly with a
given tropical polynomial and compute the homology of the $F$-sheaves on the tropical hypersurface. In other words,  the 
tropical homology of the tropical hypersurface.

\begin{example}\label{ex:conic}
We begin with a tropical curve in $\mathbb{R}^2$ which is dual to a triangulation of a square of size $1$. 
\begin{verbatim}
application "tropical";
$f = toTropicalPolynomial("max(0,x+5,y+3, x+y+9)");
$div = divisor( (projective_torus<Max>(2)), 
   rational_fct_from_affine_numerator($f));
application "fan";
@betti_usual = ();
@betti_bm = ();
for(my $i=0;$i<2;$i++){
   my $fi = $div->fcosheaf($i);
   my $si=$div->usual_chain_complex($fi);
   my $bmi=$div->borel_moore_complex($fi);
   push @betti_usual, $si->BETTI_NUMBERS;
   push @betti_bm, $bmi->BETTI_NUMBERS;
}
\end{verbatim}
gives
\begin{verbatim}
fan > print new Matrix(@betti_usual);
1 0
3 0

fan > print new Matrix(@betti_bm);
0 3
0 1
\end{verbatim}

\end{example}

\begin{example}\label{ex:K3}
As a final example, we calculate the homology of another tropical hypersurface. 
This hypersurface arises as a triangulation of the three dimensional simplex of edge length 4,
and is a tropical  K3-surface in $\mathbb{R}^3$. 

\begin{verbatim}
application "tropical";
$f = toTropicalPolynomial("max(0,x,y,z, 2*x-2, 
   2*y-2, 2*z-2, x+y-1, x+z-1, y+z-1, 3*x-6, 
   3*y-6, 3*z-6, 2*x+y-4, 2*y+x-4, 2*x+z-4, 
   2*z+x-4, 2*y+z-4, 2*z+y-4, x+y+z+1, 4*x-12, 
   4*y-12, 4*z-12, 3*x+y-9, 3*y+x-9, 3*x+z-9, 
   3*z+x-9, 3*y+z-9, 3*z+y-9, 2*x+2*y-8, 
   2*x+2*z-8, 2*y+2*z-8, 2*x+y+z-7, x+2*z+y-7, 
   2*y+z+x-7)");
$k3 = divisor((projective_torus<Max>(3)),
   rational_fct_from_affine_numerator($f));
application "fan";
@numbers = (0..2);
@cosheaves = map{$k3->fcosheaf($_)} @numbers;
@usualChainComplexes = map{$k3->usual_chain_complex($_)} 
   @cosheaves;
@bmComplexes = map{$k3->borel_moore_complex($_)} 
   @cosheaves;
@betti_usual = map{$_->BETTI_NUMBERS} 
   @usualChainComplexes;
@betti_bm = map{$_->BETTI_NUMBERS} @bmComplexes;
\end{verbatim}

We obtain the following matrices of Betti numbers:

\begin{verbatim}
fan > print new Matrix(@betti_usual);
1 0 1
3 31 0
34 0 0

fan > print new Matrix(@betti_bm);
0 0 34
0 31 3
1 0 1
\end{verbatim}

This tropical hypersurface is bigger than the polyhedral complexes we considered before.
Its $f$-vector is $(64, 96, 34)$. This can be seen from the usual chain complex of the \(F^0\)-cosheaf. 

\begin{verbatim}
fan > $usualChainComplexes[0]->print();
 3       2        1        0        -1
k^0 --> k^34 --> k^96 --> k^64 --> k^0
\end{verbatim}

In this example and  Example \ref{ex:conic}, we again observed the homological version of Poincar\'e
duality. 
\end{example}

\section{Future directions}\label{sec:future}
\subsection{Sheaves of modules}

It is also possible to compute (co)homology of cellular (co)sheaves of modules. 
For example, given a rational polyhedral complex there
are also integral versions of the $W$-sheaves and $F$-cosheaves, which are
free $\mathbb{Z}$-modules. 
However,
using the current  methods \texttt{fcosheaf} and \texttt{wsheaf} can lead to
incorrect results over $\ZZ$. Still, the ranks of the torsion and the free part
of the (co)homology will be correct in these cases.

The problem with using the current implementation to compute integral versions of
the (co)homology of the integral versions of these (co)sheaves is that the
property \texttt{BASES} does not necessarily consist of a lattice basis of the free
$\ZZ$-module for each face. In addition, the matrices in \texttt{BLOCKS} may not 
accurately encode the $\ZZ$-linear maps. 

If one properly chooses $\ZZ$-bases for \texttt{BASES} and 
defines \texttt{BLOCKS} manually  with the correct
maps over $\ZZ$ when creating a (co)sheaf, then the current rules for computing the
cellular (co)homology will compute the correct $\ZZ$-homology. 

We plan to adapt \texttt{fcosheaf} and \texttt{wsheaf} to give the correct
results over $\ZZ$ after switching to \polymake 's internal chain complex
object. This has  recently been pushed to the \polymake repository by Olivia
R\"ohrig.

\subsection{Tropical compactifications and projective hypersurfaces}\label{sec:compact}
How to implement compact tropical varieties is part of an ongoing discussion
inside the \polymake developer team. One possibility is to save one affine
tropical variety per chart of the tropical projective space. For many cases
this would result in a drastic increase of resource usage. Thus, one may want
to restrict to certain classes of tropical varieties with nice
compactifications. 

A solution to this problem is necessary  in order to use our extension in order to 
give answers for example  Problem 10 on Surfaces in \cite{Sturmfels}. Upon having an 
implementation of compact tropical varieties, one could for example combine 
our package and the 
approach to tropical Enriques surfaces in \cite{Enriques} to determine the Hodge numbers 
and solve Problem 10.

\subsection{Implementing other cellular (co)sheaves}
There are  many other imaginable cellular (co)sheaves to consider on a
polyhedral complex, including  (co)sheaves arising from common
(co)sheaf operations, such as restrictions, pullbacks, and Verdier duals. 

Cellular sheaves on polyhedral fans have appeared in the work of Brion (see
\cite{Brion}). 
There, given a polyhedral fan in $\mathbb{R}^n$, one associates to a
face $\sigma$,  the vector space $\mathbb{R}^n / Lin(\sigma)$, where
$Lin(\sigma)$ denotes the linear span of the $\sigma$.  These vector spaces
come equipped with  natural maps between them when there is an inclusion of
faces. One can also take $p^{\text{th}}$ exterior powers of these vector
spaces, as well as generalise this definition beyond polyhedral fans to get a
collection of sheaves.

The cohomology of these cellular sheaves is related to the motion spaces in
discrete dynamical geometry (see \cite{Whiteley}). Following the descriptions
in that paper, the  $W$-sheaves come up in aspects of rigidity and the
$F$-cosheaves of skeleta of polyhedral complexes are related to stress spaces.

\subsection{Applied topology}

Celluar (co)sheaves have also appeared often in the field of applied topology and 
topological data analysis, notably  in the study of sensor networks, network coding, and 
persistent homology (see \cite{Curry},  \cite{Ghrist}). 
Although the (co)sheaves appearing in these contexts are often a part of the input data of the model under consideration and do not 
have a simple recipe coming from the geometry of the underlying topological space like in the case of tropical homology. 
However, our  extension allows for the construction of a sheaf from scratch. 
Another generalisation to be considered in the future, is that the underlying topological spaces appearing in this context are not necessarily polyhedral complexes in 
$\mathbb{R}^n$. The current extension capabilities for (co)sheaves on polyhedral complexes could be extended to more general topological spaces 
using the \polymake application \texttt{topaz}. 

We would also like to point out the current efforts underway by Olivia R\"ohrig  to implement persistent homology in \polymake.

\section*{Acknowledgements}
  
This package was developed while all three authors were visiting the Fields
Institute for Research in Mathematical Sciences in Toronto during the thematic
program on Combinatorial Algebraic Geometry. We are very grateful to the
organisers of the semester and the institute for their hospitality. 

We would like to thank  Greg Smith, Bernd Sturmfels and  five anonymous referees for their careful attention to an 
earlier version of this manuscript.

\footnotesize
\noindent {\bf Authors' addresses:}

\smallskip

\noindent 
Lars Kastner,  Freie Universit\"at Berlin, Arnimallee 3, 14195 Berlin, Germany, { \tt k.l@fu-berlin.de}

\smallskip

\noindent Kristin Shaw Technische Universit\"at Berlin, Institut f\"ur Mathematik, Stra\ss e des 17. Juni 136, 10623 Berlin, Germany, {\tt shaw@math.tu-berlin.de}

\smallskip

\noindent  Anna-Lena Winz Freie Universit\"at Berlin, Arnimallee 3, 14195 Berlin, Germany, {\tt anna-lena.winz@fu-berlin.de}

\end{document}